\documentclass{birkjour}

\usepackage[all,cmtip]{xy}
 \newtheorem{thm}{Theorem}[section]
 \newtheorem{cor}[thm]{Corollary}
 
 \newtheorem{prop}[thm]{Proposition}
 \theoremstyle{definition}
 \newtheorem{defn}[thm]{Definition}
 \theoremstyle{remark}
 \newtheorem{rem}[thm]{Remark}
 \newtheorem*{ex}{Example}
 \numberwithin{equation}{section}
 
 \newcommand*{\bR}{\ensuremath{\mathbb{R}}}
 \newcommand*{\R}{\ensuremath{\mathbb{R}}}
 
 \DeclareMathOperator{\diam}{diam}
 \DeclareMathOperator{\Modd}{Mod}

\DeclareMathOperator{\modulus}{Mod}
\DeclareMathOperator{\loc}{loc}

\def\Xint#1{\mathchoice
	{\XXint\displaystyle\textstyle{#1}}%
	{\XXint\textstyle\scriptstyle{#1}}%
	{\XXint\scriptstyle\scriptscriptstyle{#1}}%
	{\XXint\scriptscriptstyle\scriptscriptstyle{#1}}%
	\!\int}
\def\XXint#1#2#3{{\setbox0=\hbox{$#1{#2#3}{\int}$}
		\vcenter{\hbox{$#2#3$}}\kern-.5\wd0}}

\def\dashint{\Xint-}

\begin{document}

%
%
%
%
%
%
%
%
%

\title[QR mappings in metric spaces]
 {The theory of quasiregular mappings in metric spaces: progress and challenges}

\author[C.-Y. Guo]{Chang-Yu Guo}

\address{%
Department of Mathematics\\
University of Fribourg\\
Chemin du Musee 23\\
CH-1700 Fribourg, Switzerland}

\email{changyu.guo@unifr.ch}

\thanks{This work was completed with the support of the Swiss National Science Foundation (No. 153599).}



\date{\today}
\dedicatory{Dedicated to Seppo Rickman and Jussi V\"ais\"al\"a on the occasion of their $80^\text{th}$ birthdays }

\begin{abstract}
We survey the recent developments in the theory of quasiregular mappings in metric spaces. In particular, we study the geometric porosity of the branch set of quasiregular mappings in general metric measure spaces, and then, introduce the various natural definitions of quasiregular mappings in general metric measure spaces, and give conditions under which they are quantitatively equivalent.
\end{abstract}

\maketitle
\section{Introduction}

A continuous mapping $f\colon X\to Y$ between topological spaces is said to be a \textit{branched covering} if $f$ is \textit{discrete}, \textit{open} and \textit{of locally bounded multiplicity}. Recall that $f$ is open if it maps each open set in $X$ to an open set $f(X)$ in $Y$; $f$ is discrete if for each $y\in Y$ the preimage $f^{-1}(y)$ is a discrete subset of $X$; $f$ has locally bounded multiplicity if for each $x\in X$, there exists a neighborhood $U_x$ of $x$ and a positive constant $M_x$ such that $N(f,U_x)\leq M_x<\infty$, where $N(f,U_x):=\sup_{y\in Y}\sharp\{f^{-1}(y)\cap U_x\}$ is the multiplicity of $f$ in $U_x$. The latter assumption can be dropped off if $X$ and $Y$ have certain manifold structure~\cite{hr02}. Nonconstant holomorphic functions between connected Riemann surfaces are typical examples of branched coverings.

For a branched covering $f\colon X\to Y$ between two metric spaces, $x\in X$ and $r>0$, set
\begin{equation*}
H_f(x,r)=\frac{L_f(x,r)}{l_f(x,r)},
\end{equation*}
where
\begin{equation*}
L_f(x,r):=\sup_{y\in Y}{\{d(f(x),f(y)):d(x,y)= r\}},
\end{equation*}
and
\begin{equation*}
l_f(x,r):=\inf_{y\in Y}{\{d(f(x),f(y)):d(x,y)= r\}}.
\end{equation*}
Then the \textit{linear dilatation function} of $f$ at $x$ is defined pointwise by 
\begin{equation*}
H_f(x)=\limsup_{r\to0}H_f(x,r).
\end{equation*}

\begin{defn}\label{def:metric quasiregular map}
	A branched covering $f\colon X\to Y$ between two metric measure spaces\footnote{A \textit{metric measure space} is defined to be a triple $(X,d,\mu)$, where $(X,d)$ is a \emph{separable} metric space and $\mu$ is a \emph{nontrivial locally finite Borel regular measure} on $X$.} is termed \textit{metrically $H$-quasiregular} if the linear dilatation function $H_f$ is finite everywhere and essentially bounded from above by $H$. 
\end{defn}

If $f\colon X\to Y$, in Definition~\ref{def:metric quasiregular map}, is additionally assumed to be a homeomorphism, then $f$ is called \textit{metrically $H$-quasiconformal}. We will call $f$ a metrically quasiregular or quasiconformal mapping if it is metrically $H$-quasiregular or $H$-quasiconformal for some $H\in [1,\infty)$.

Quasiregular mappings were first introduced by Reshetnyak in 1966~\cite{r66}, where he actually used the (equivalent) analytic formulation of quasiregularity. The analytic foundation of the theory of quasiregular mappings were laid after a sequence of his papers from 1966 to 1969. A deep fact he discovered is that analytic quasiregular mappings are branched coverings~\cite{re89}. The whole theory of quasiregular mapping were significantly advanced in a sequence of papers from the Finnish school in the late 1960s~\cite{mrv69,mrv70,mrv71}. See also the beautiful paper~\cite{m14} for a nice survey on the development of the field.

For quasiconformal mappings, there is also a geometric definition, which makes uses of the modoulus of curve families. One of the most fundamental result in the theory of quasiconformal mappings, due to the deep works of Gehring, V\"ais\"al\"a and many others, is that all the three different definitions of quasiconformality are quantitatively equivalent. According to a later remarkable result of Heinonen and Koskela~\cite{hk95}, one can even relax limsup to liminf in the metric definition of quasiconformality. The similar result holds in the theory of quasiregular mappings as well, but the proofs are much more involved~\cite{r93}. 

Due the numerious successful applications of the theory quasiconformal mappings in geometric group theory (for instance~\cite{m73}), the fundation of the theory of metrically quasiconformal mappings has been laid~\cite{hk98} in the general framework of metric measure spaces with controlled geometry.

All of the three definitions of quasiconformality can be generalized in a natural way to the setting of metric measure spaces. A remarkable fact, after a number of seminal works~\cite{hkst01,t98,t00,t01,bkr07,w12proc,w14}, is that the quantitative equivalences of quasiconformality extends to a large class of metric spaces.

The recent advances in analysis on metric spaces~\cite{hr02,hs02,hk11} promotes a general theory of quasiregular mappings in the setting of metric spaces, whereas a complete understanding of the various definitions of quasiregularity plays an important role. In this survey, we examine
the minimal assumptions on metric spaces $X$ and $Y$ and a branched covering $f\colon X\to Y$  guarantee $f$ being quasiregular according to different definitions.

\section{Preliminaries on metric spaces}

A main theme in analysis on metric spaces is that the infinitesimal
structure of a metric space can be understood via the curves that it contains.
The reason behind this is that we can integrate Borel functions along rectifiable curves and do certain non-smooth calculus akin to the Euclidean spaces.

Let $(X,d)$ be a metric space. A curve (or path) in $X$ is a continuous map $\gamma\colon I\to X$, where $I\subset \bR$ is an interval. We call $\gamma$ compact, open, or half-open, depending on the type of the interval $I$.

Given a compact curve $\gamma\colon [a,b]\to X$, we define the \textit{variation function} $v_\gamma\colon [a,b]\to [a,b]\to \infty$ by
\begin{align*}
v_\gamma(s)=\sup_{a\leq a_1\leq b_1\leq \cdots\leq a_n\leq b_n\leq s}\sum_{i=1}^n d(\gamma(b_i),\gamma(a_i)).
\end{align*}
The length $l(\gamma)$ of $\gamma$ is defined to be the variation $v_\gamma(b)$ at the end point $b$ of the parametrizing interval $[a,b]$. If $\gamma$ is not compact, its length $l(\gamma)$ is defined to be the supremum of the lengths of the compact subcurves of $\gamma$. 

A curve is said to be \textit{rectifiable} if its length $l(\gamma)$ is finite, and \textit{locally rectifiable} if each of its compact subcurves is rectifiable. For any rectifiable curve $\gamma$ there are its associated length function $s_\gamma\colon I\to [0,l(\gamma)]$ and a unique 1-Lipschitz map $\gamma_s\colon [0,l(\gamma)]\to X$ such that $\gamma=\gamma_s\circ s_\gamma$. The curve $\gamma_s$ is the \textit{arc length parametrization} of $\gamma$.

When $\gamma$ is rectifiable, and parametrized by arclength on the interval $[a,b]$, the integral of a Borel function $\rho\colon X\to [0,\infty]$ along  $\gamma$ is
\[ \int_\gamma \rho \,ds = \int_0^{l(\gamma)} \rho(\gamma_s(t))\,dt\text{.}\]
Similarly, the line integral of a  Borel function $\rho\colon X\to [0,\infty]$ over a locally rectifiable curve $\gamma$ is defined to be the supremum of the integral of $\rho$ over all compact subcurves of $\gamma$.

A curve $\gamma$ is \textit{absolutely continuous} if $v_\gamma$ is absolutely continuous. Via the chain rule, we then have
\begin{align}\label{eq:formula for ABS curve line integral}
\int_{\gamma}\rho ds=\int_a^b\rho(\gamma(t))v_\gamma'(t)dt.
\end{align}

Let $X=(X,d,\mu)$ be a metric measure space. Let $\Gamma$ a family of curves in $X$.  A Borel function $\rho\colon X\rightarrow [0,\infty]$ is \textit{admissible} for $\Gamma$ if for every locally rectifiable curves $\gamma\in \Gamma$,
\begin{equation}\label{admissibility}
\int_\gamma \rho\,ds\geq 1\text{.}
\end{equation}
The \textit{$p$-modulus} of $\Gamma$, $p\geq 1$, is defined as
\begin{equation*}
\modulus_p(\Gamma) = \inf_{\rho} \left\{ \int_X \rho^p\,d\mu:\text{$\rho$ is admissible for $\Gamma$} \right\}.
\end{equation*}
A family of curves is called \emph{$p$-exceptional} if it has $p$-modulus zero. We say that a property of curves holds for \emph{$p$-almost every curve} if the collection of curves for which the property fails to hold is $p$-exceptional.

Let $X=(X,d,\mu)$ be a metric measure space and $Z=(Z,d_Z)$ be a metric space.

A Borel function $g\colon X\rightarrow [0,\infty]$ is called an \textit{upper gradient} for a map $f\colon X\to Z$ if for every rectifiable curve $\gamma\colon [a,b]\to X$, we have the inequality
\begin{equation}\label{ugdefeq}
\int_\gamma g\,ds\geq d_Z(f(\gamma(b)),f(\gamma(a)))\text{.}
\end{equation}
If inequality \eqref{ugdefeq} holds for $p$-almost every compact curve, then $g$ is called a \textit{$p$-weak upper gradient} for $f$.  When the exponent $p$ is clear, we omit it.

A $p$-weak upper gradient $g$ of $f$ is \textit{minimal} if for every $p$-weak upper gradient $\tilde{g}$ of $f$, $\tilde{g}\geq g$ $\mu$-almost everywhere.  If $f$ has an upper gradient in $L^p_{\loc}(X)$, then $f$ has a unique (up to sets of $\mu$-measure zero) minimal $p$-weak upper gradient.  We denote the minimal upper gradient by $|\nabla f|$.

Real-valued Sobolev spaces based on upper gradients were used to great
success~\cite{c99}  and explored in-depth in~\cite{s00}. They have been extended to the
metric-valued setting in~\cite{hkst01,hkst15}. There are several equivalent ways to define the Sobolev spaces of mappings between metric measure space and a simple definition is as follows. Let $f\colon X\to Y$ be a continuous map. Then $f$ belongs to the Sobolev space $N^{1,p}_{\loc}(X,Y)$, $1\leq p<\infty$, if for each relatively compact open subset $U\subset X$, the map $f$ has an upper gradient $g\in L^p(U)$ in $U$, and
there is a point $x_0\in U$ such that $u(x):=d_Y(u(x),x_0)\in L^p(U)$. 

A Borel regular measure $\mu$ on a metric space $(X,d)$ is called a \textit{doubling measure} if every ball in $X$ has positive and finite measure and there exists a constant $C_\mu\geq 1$ such that
	\begin{equation}\label{eq:doubling measure}
	\mu(B(x,2r))\leq C_\mu \mu(B(x,r))
	\end{equation}
for each $x\in X$ and $r>0$. We call the triple $(X,d,\mu)$ a doubling metric measure space if $\mu$ is a doubling measure on $X$. We call $(X,d,\mu)$ an \textit{Ahlfors $Q$-regular space}, $1\leq Q<\infty$, if there exists a constant $C\geq 1$ such that
	\begin{equation}\label{eq:Ahlfors regular measure}
	C^{-1}r^Q\leq \mu(B(x,r))\leq Cr^Q
	\end{equation}
for all balls $B(x,r)\subset X$ of radius $r<\diam X$.

We say that a metric measure space $(X,d,\mu)$ admits a \textit{(1,p)-Poincar\'e inequality} if there exist constants $C\geq 1$ and $\tau\geq 1$ such that
\begin{equation}\label{eq:Poin ineq}
\dashint_B|u-u_B|d\mu\leq C\diam(B)\Big(\dashint_{\tau B}g^pd\mu\Big)^{1/p}
\end{equation}
for all open balls $B$ in $X$, for every function $u:X\to\bR$ that is integrable on balls and for every upper gradient $g$ of $u$ in $X$.

The $(1,p)$-Poincar\'e inequality can be thought of as a requirement that a
space contains ``many" curves, in terms of the $p$-modulus of curves in the
space. For more information on the Poincar\'e inequalities, see~\cite{hk00,hkst15}.

A metric measure space $(X,d,\mu)$ is said to have \textit{$Q$-bounded geometry} if it is Ahlfors $Q$-regular and supports a $(1,Q)$-Poincar\'e inequality. See~\cite{hk98} for more information on metric spaces of bounded geometry.

Throughout the paper, we assume that $X$ and $Y$ are \textit{locally compact}, \textit{complete}, \textit{connected}, \textit{locally connected} metric spaces

\section{Size of the branch set}\label{sec:size of branch set}

The main obstacle in establishing the theory of quasiregular mappings in general metric spaces lies in the \textit{branch set} $\mathcal{B}_f$, i.e., the set of points in $X$ where $f\colon X\to Y$ fails to be a local homeomorphism. The difficulty is somehow hidden in the Euclidean planar case, as the celebrated \textit{Stoilow factorization theorem} asserts that a quasiregular mapping $f\colon \Omega\to \R^2$
admits a factorization $f=\varphi\circ g$, where $g\colon \Omega\to g(\Omega)$ is quasiconformal and $\varphi\colon g(\Omega)\to \R^2$ is analytic. This factorization, together with relatively complete understanding of the structure of analytic functions in the plane, connects quasiregular and quasiconformal mappings strongly. In particular, the branch set $\mathcal{B}_f$ of a quasiregular mapping $f\colon \Omega\to \R^2$ is discrete.

In higher dimensions or more general metric measure spaces, the branch set of a quasiregular mappings can be very wild, for instance, it might contain many wild Cantor sets, such as the Antoine's necklace~\cite{hr98}, of classical geometric topology. This makes the homeomorphic theory and the non-homeomorphic theory substantially different. Indeed,  the most delicate part of establishing the theory of quasiregular mappings in various settings as mentioned above is to show that the branch set and its image have null measure. For a survey on the topological property of the branch set as well as open problems in this direction, see~\cite{h02}. 

Regarding the Hausdorff dimension of $\mathcal{B}_f$ and its image $f\big(\mathcal{B}_f\big)$ in the Euclidean setting, a well-known result of Gehring and V\"ais\"al\"a~\cite{gv72} says that for each $n\geq 3$ and each pair of numbers $\alpha,\beta\in [n-2,n)$, there exists a quasiregular mapping $f:\bR^n\to \bR^n$ such that 
\begin{align*}
\dim_{\mathcal{H}}\mathcal{B}_f=\alpha\quad \text{and}\quad \dim_{\mathcal{H}}f(\mathcal{B}_f)=\beta.
\end{align*}
On the other hand, by the result of Sarvas~\cite{s75}, for a non-constant $H$-quasiregular mapping $f\colon \Omega\to \bR^n$, $n\geq 2$, between Euclidean domains, 
\begin{align}\label{eq:Sarvas}
\dim_{\mathcal{H}}f(\mathcal{B}_f)\leq n-\eta
\end{align}
for some constant $\eta=\eta(n,H)>0$. Yet another well-known result of Bonk and Heinonen~\cite{bh04} says that for a non-constant $H$-quasiregular mapping $f\colon \Omega\to \bR^n$, $n\geq 2$, between Euclidean domains, 
\begin{align}\label{eq:Bonk-Heinonen}
\dim_{\mathcal{H}}\mathcal{B}_f\leq n-\eta
\end{align}
for some constant $\eta=\eta(n,H)>0$.

The recent development in analysis on metric spaces~\cite{hr02,hs02,hk11} promotes a general theory of quasiregular mappings beyond the Riemannian spaces, whereas a deeper understanding of the branch set of a quasiregular mapping is rather necessary. 

Assume that $X$ and $Y$ are doubling metric spaces, which are also topological $n$-manifolds, that $X$ is linearly locally contractible\footnote{A metric space is \textit{linearly locally contractible} if each ball of raidus $r$ contracts in a ball with the same center and radius $\lambda r$ for some $\lambda\geq 1$.}, and that $Y$ has bounded turning\footnote{A metric space $X$ has \textit{$c$-bounded turning} if every pair of points $x_1,x_2\in X$ can be joined by a continuum $E\subset X$ such that $\diam E\leq cd(x_1,x_2)$.}. A special case of~\cite[Theorem 1.1]{gw15} reads as follows.

\begin{thm}\label{thm:Theorem 1.1}
	If $H_f(x)\leq H$ for every $x\in X$, then $\mathcal{B}_f$ and $f(\mathcal{B}_f)$ are countably $\delta$-porous, quantitatively. Moreover, the porosity constant can be explicitly calculated.
\end{thm}

Recall that a set $E\subset X$ is said to be $\alpha$-\textit{porous} if for each $x\in E$, 
\begin{align*}
\liminf_{r\to 0}r^{-1}\sup\big\{\rho:B(z,\rho)\subset B(x,r)\backslash E\big\}\geq \alpha.
\end{align*}
A subset $E$ of $X$ is called  \textit{countablely ($\sigma$-)porous} if it is a countable union of ($\sigma$-)porous subsets of $X$.

Since porous sets have Hausdorff dimension strictly less than $Q$, quantitative, in an Ahflors $Q$-regular space, we have the following immediate corollary.
\begin{cor}\label{coro:Corollary 1.3}
	If $X$ and $Y$ are Ahlfors $Q$-regular, and $H_f(x)<\infty$ for all $x\in X$, then $\mathcal{H}^Q(\mathcal{B}_f)=\mathcal{H}^Q(f(\mathcal{B}_f))=0$. Moreover, if $H_f(x)\leq H$  for all $x\in X$, then 
	\begin{align*}
	\max\big\{\dim_{\mathcal{H}}(\mathcal{B}_f),\dim_{\mathcal{H}}(f(\mathcal{B}_f))\big\}\leq Q-\eta<Q,
	\end{align*}
	where $\eta$ depends only on $H$ and the data of $X$ and $Y$. Moreover, $\eta$ can be explicitly calculated.
\end{cor}

The assumptions in Theorem~\ref{thm:Theorem 1.1} or Corollary~\ref{coro:Corollary 1.3} is rather sharp as the following example indicates.

\begin{ex}[\cite{gw16}]
	For each $n\geq 3$, there exist an Ahlfors $n$-regular metric space $X$ that is homeomorphic to $\bR^n$ and supports a $(1,1)$-Poincar\'e inequality, and a 1-quasiregular mapping (indeed 1-BLD) $f\colon X\to \bR^n$, such that $\min\big\{\mathcal{H}^n(\mathcal{B}_f),\mathcal{H}^n(f(\mathcal{B}_f))\big\}>0.$  
\end{ex}


\section{The pullback factorization}\label{sec:the pullback factorization}
In this section, we briefly explore the pullback factorization introduced (due to M. Williams) in~\cite{gw16}. One feature of the Stoilow facterization is that we can factorize a quasiregular mapping $f$ into the composition of an ``analytically nice" (high regularity) mapping $\pi$ with a ``topologically nice" (homeomorphism) mapping $g$. 

Let $f\colon X\to Y$ be a branched covering.

\subsection{The pullback metric}\label{subsec:The pullback metric}
The ``\textit{pullback metric}" $f^*d_Y\colon X\times X\to [0,\infty)$ is defined as follows:
\begin{align}\label{eq:def for pullback metric}
f^*d_Y(x_1,x_2)=\inf_{x_1,x_2\in \alpha}\diam\big(f(\alpha)\big),
\end{align}
where the infimum is taken over all continua $\alpha$ joining $x_1$ and $x_2$ in $X$.

It is immediate from the definition that $f^*d_Y$ satisfies the triangle inequality. Moreover, the connectivity assumption on $X$ guarantees that it is finite. Thus the discreteness of $f$ implies that $f^*d_Y$ is indeed a genuine metric. We denote by $X^f$ or $f^*Y$ the metric space $(X,f^*d_Y)$.


\subsection{Canonical factorization}\label{subsec:Canonical factorization}
In what follows, let $g\colon X\to X^f$ be the identity map, and let $\pi\colon X^f\to Y$ satisfy $\pi\circ g=f$, so that on the level of sets, $f=\pi$. We refer this canonical factorization as the \textit{pullback factorization for} $f$. 

\begin{figure}[h]
	\[
	\xymatrix{
		X \ar[r]^{g} \ar[dr]_{f} & X^f \ar[d]^{\pi} \\
		& Y
	}
	\]
	\caption{The canonical pullback factorization}\label{Fig:pullback factorization}
\end{figure}
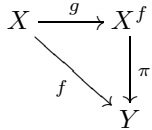

Under appropriate conditions, the metric space $X^f$ and the branched covering $\pi$ are well behaved. In order to make precise these nice behaviors, we need to recall some of the nice mapping classes between metric spaces.

\begin{defn}\label{def:BLD mapping}
	A branched covering $f\colon X\to Y$ between two metric spaces is said to be an \textit{$L$-BLD}, or a \textit{mapping of $L$-bounded length distortion}, $L\geq 1$, if 
	\begin{align*}
	L^{-1}l(\alpha)\leq l(f\circ \alpha)\leq Ll(\alpha)
	\end{align*}
	for all non-constant curves $\alpha$ in $X$.
\end{defn}
The definition of BLD mappings is clearly only interesting if the metric spaces $X$ and $Y$ have a reasonable supply of rectifiable curves and so the most natural setting in which we study such mappings is that of quasiconvex metric spaces\footnote{A metric space $X$ is \textit{$c$-quasiconvex} if every
	pair of points $x_1,x_2\in X$ may be joined with a curve $\gamma$ of length $l(\gamma)\leq cd(x_1,x_2)$. }.


When $X$ and $Y$ have bounded turning, the natural branched analog of a bi-Lipschitz homeomorphism is what we call a \textit{mapping of bounded diameter distortion}, which is defined in analogy with BLD mappings:
\begin{defn}\label{def:BDD mapping}
	A branched covering $f\colon X\to Y$ between two metric spaces is said to be an \textit{$L$-BDD}, or a \textit{mapping of $L$-bounded diameter distortion}, $L\geq 1$, if 
	\begin{align*}
	L^{-1}\diam(\alpha)\leq \diam(f\circ \alpha)\leq L\diam(\alpha)
	\end{align*}
	for all non-constant curves $\alpha$ in $X$.
\end{defn}

It follows directly from the definition of arc-length that an $L$-BDD mapping is $L$-BLD as well. 

As we will see in a moment, the metric space $X^f$ retains the original topology of $X$ and is often rather well-behaved: having 1-bounded turning and inheriting many metric and geometric properties from $Y$. Moreover, the branched covering $\pi$ from the pullback factorization is easily seen to be 1-Lipschitz and 1-BDD (and a fortiori 1-BLD).

\subsection{Fine properties of the pullback metric}
In this section, we fix a branched covering $\pi\colon X\to Y$ between two metric spaces. However, we will typically consider the space $X^\pi$ (or $\pi^*Y$) by endowing the set $X$ with the pullback metric $\pi^*d_Y$.

For each metric $d$ on a topological space $X$, the length metric $l_d(z_1,z_2)$ is given by infimizing the lengths of all paths joining $z_1$ and $z_2$. For a metric space $X=(X,d_X)$, we denote by $X^l$ the length space $(X,l_{d_X})$.


For simplicity, we will formulate many of our basic results for the case that $\pi$ is \textit{proper} (i.e.~the preimage of each compact set is compact) and surjective, and $N=N(\pi,X)<\infty$. Up to the end of this section, we additionally assume the metric space $Y$ is \textit{proper}, i.e. closed bounded balls are compact. 

Recall that a continuous mapping $f\colon X\to Y$ between two metric spaces is said to be \textit{$c$-co-Lipschitz} if for all $x\in X$ and $r>0$, 
\begin{align*}
B(f(x),r)\subset f(B(x,cr)).
\end{align*}

The basic properties of the pullback metric are formulated in the next proposition. Recall that for a continuous mapping $\pi\colon X\to Y$, $U(x,\pi,r)$ denotes the $x$-component of $\pi^{-1}(B(\pi(x),r))$.

\begin{prop}[Section 4.3, \cite{gw16}]\label{prop:basic properties of pullback metric}
	
\begin{itemize}
	\item[1).] The metric space $X^\pi$ is a proper metric space, homeomorphic to $X$ via the indentity mapping $g$. Open and closed balls in $X^\pi$ are connected and $X^\pi$ have 1-bounded turning. The projection mapping $\pi\colon X^\pi\to Y$ is 1-Lipschitz, 1-BDD and for each $z\in X^\pi$,
	\begin{align}\label{eq:pullback property 1}
	B(z,r)\subset U(z,\pi,r)\subset B(z,2r).
	\end{align}
	
	\item[2).] 	If $Y$ has $c$-bounded turning, then $\pi$ is $c$-co-Lipschitz and is locally $c$-bi-Lipschitz on each set $X_k:=\{z\in X^\pi:N(z,\pi,X)=k\}$. If, additionally, $Y$ is Ahlfors $Q$-regular with constant $c_2$, then $X^\pi$ is Ahlfors $Q$-regular with constant $c^Qc_2N$. Moreover, for each $k=1,\dots,N$, and at each Lebesgue point $z$ of $X_k$, we may take the pointwise constant of $Q$-regularity to be $c^Qc_2$.
	
	\item[3).] If $Y$ is $c$-LLC-2 and $X$ has no local cut points\footnote{A point $x$ in a metric space $X$ is called a \textit{local cut point} if $U\backslash \{x\}$ is disconnected for some neighborhood $U$ of $x$}, then $X^\pi$ is locally $2c$-LLC.
	
	\item[4).] If $l_{d_Y}$ and $d_Y$ induces the same topology on $Y$, then the metrics $\pi^*d_Y$, $l_{\pi^*d_Y}$, $\pi^*l_{d_Y}$, and $l_{\pi^*l_{d_Y}}$ induces the same topology on $X$. The length of a curve in $X^\pi$ is the same with respect to any of these four metrics, and moreover,
	\begin{align*}
	\pi^*d_Y\leq \pi^*l_{d_Y}\leq l_{\pi^*d_Y}=l_{\pi^*l_{d_Y}}\leq (2N-1)\pi^*l_{d_Y}.
	\end{align*}
	In particular, the metric space $\pi^*(Y^l)$ is $(2N-1)$-quasiconvex and if $Y$ is $c$-quasiconvex, then $X^\pi$ is $(2N-1)c$-quasiconvex. 
\end{itemize}

\end{prop}

\subsection{Fine properties of the pullback factorization}

Fix any branched covering $f\colon X\to Y$ and the pullback factorization $f=\pi\circ g$, where
$g\colon X\to X^f$ is the identity mapping and $\pi\colon X^f\to Y$ is the branched covering given by $\pi=f$. 

By Proposition~\ref{prop:basic properties of pullback metric}, $\pi$ is a 1-BDD mapping and thus we have factorized $f$ into a composition of a homeomorphism and a 1-BDD ``projection". Note that while on the level of sets, we are factoring out the identity, the mapping $g$ will typically not be an isometry. But we have seen already, the projection mapping $\pi$ can be thought of as being as close to an isometry as possible. Thus, philosophically, we have factored $f$ into a geometric equivalence composed with a topological one.

Furthermore, as a result of the fact that $\pi$ is 1-BDD, $f$ and $g$ share many geometric properties.
\begin{prop}[Section 4.4, \cite{gw16}]\label{prop:nice properties of pullback factorization I}
\begin{itemize}
\item[1).] The mapping $f$ is $L$-BDD/BLD if and only if $g$ is $L$-BDD/BLD. 

\item[2).] Suppose that $Y$ is $c$-quasiconvex, and that $f$ is $L$-Lipschitz and $L$-BLD with $N=N(f,X)<\infty$. Then $f$ is $cNL$-BDD. 
\end{itemize}
\end{prop}

Thus much of the theory of BLD and BDD mappings reduces to the study of the pullback metric. Similarly, it turns out that under various definitions and for many different levels of generality, $g$ is quasiconformal if and only if $f$ is quasiregular. Thus one obtains a canonical factorization of a quasiregular mapping into a composition of a quasiconformal mapping
with a 1-BDD mapping. This is particularly useful in extending the theory of quasiregular mappings to the metric setting, as the quasiconformal theory has at present advanced much further than its branched counterpart in this generality. 

\section{Quasiregular mappings between metric spaces}

Fix two metric measure spaces $X=(X,d_X,\mu)$ and $Y=(Y,d_Y,\nu)$. Assume that $\tilde{\Omega}\subset X$ is a domain, and that $f\colon \tilde{\Omega}\to \Omega\subset Y$ is an \textit{onto} branched covering. 

In what follows, we fix the pullback factorization $f=\pi\circ g$ as in Section~\ref{sec:the pullback factorization}, where $g\colon \tilde{\Omega}\to \tilde{\Omega}^f$ is the identity mapping and $\pi\colon \tilde{\Omega}^f\to \Omega$ is the 1-BDD projection. We equip $\tilde{\Omega}^f$ with the Borel regular measure $\lambda=\pi^*\nu=g_*f^*\nu$ and simply write $\tilde{\Omega}^f$ for the metric measure space $(\tilde{\Omega}^f,f^*d_Y,\lambda)$.

\begin{figure}[h]
	\[
	\xymatrix{
		(\tilde{\Omega},\mu) \ar[r]^{g} \ar[dr]_{f} & (\tilde{\Omega}^f,\lambda) \ar[d]^{\pi} \\
		& (\Omega,\nu)
	}
	\]
	\caption{The new pullback factorization}
\end{figure}


\subsection{Definitions of quasireguarity in general metric measure spaces}\label{subsec:Definitions of quasireguarity in general metric measure spaces}
In this section, we introduce the different definitions of quasiregularity in general metric measure spaces. 

Set $h_f(x)=\liminf_{r\to 0}H_f(x,r)$.
\begin{defn}[Weak metrically quasiregular mappings]\label{def:weak metric qr}
	A branched covering $f\colon \tilde{\Omega}\to \Omega$ is said to be \textit{weakly metrically $H$-quasiregular} if it satisfies 
	\begin{itemize}
		\item[i).] $h_f(x)<\infty$ for all $x\in \tilde{\Omega}$;
		\item[ii).] $h_f(x)\leq H$ for $\mu$-almost every $x\in \tilde{\Omega}$. 
	\end{itemize}
\end{defn}

For a branched covering $f\colon X\to Y$, the volume Jacobian is defined by
$$J_f:=\frac{d(f^*\nu)}{d\mu},$$
where the pullback measure $f^*\nu$ on $X$ is given by 
$$f^*\nu(A)=\int_YN(y,f,A)d\nu(y).$$
\begin{defn}[Analytically quasiregular mappings]\label{def:analytic def for qr}
	A branched covering $f\colon \tilde{\Omega}\to \Omega$ is said to be \textit{analytically $K$-quasiregular} if $f\in N^{1,Q}_{loc}(\tilde{\Omega},\Omega)$ and 
	\begin{align*}
	|\nabla f|(x)^Q\leq KJ_f(x)
	\end{align*}
	for $\mu$-a.e. $x\in \tilde{\Omega}$. 
\end{defn}

The geometric definition requires some modulus inequalities between curve families.

\begin{defn}[Geometrically quasiregular mappings]\label{def:geometric def for qr}
	A branch covering $f\colon \tilde{\Omega}\to \Omega$ is said to be \textit{geometrically $K$-quasiregular} if it satisfies the \textit{$K_O$-inequailty}, i.e., for each open set $\tilde{\Omega}_0\subset \tilde{\Omega}$ and  each path family $\Gamma$ in $\tilde{\Omega}_0\subset \tilde{\Omega}$, if $\rho$ is a test function for $f(\Gamma)$, then
	\begin{align*}
	\Modd_Q(\Gamma)\leq K\int_{\Omega}N(y,f,\tilde{\Omega}_0)\rho^Q(y)d\nu(y).
	\end{align*}
\end{defn}

We will refer to the metric definition $(M)$, the weak metric definition $(m)$, the analytic definition $(A)$, and the geometric definition $(G)$ as elements of the \textit{forward definitions}.

Next, we introduce the elements from the \textit{inverse definitions}: the inverse metric definition $(M^*)$, the inverse weak metric definition $(m^*)$, the inverse analytic definition $(A^*)$, and the inverse geometric definition $(G^*)$.

For each $x\in X$, 
\begin{align*}
H_f^*(x,s)=\frac{L_f^*(x,s)}{l_f^*(x,s)},
\end{align*}
where
\begin{align*}
L_f^*(x,s)=\sup_{z\in \partial U(x,f,s)}d(x,z)\quad \text{and}\quad l_f^*(x,s)=\inf_{z\in \partial U(x,f,s)}d(x,z).
\end{align*}
The \textit{inverse linear dilatation function} of $f$ at $x$ is defined pointwise by
\begin{align*}
H_f^*(x)=\limsup_{s\to 0}H_f^*(x,s).
\end{align*}
Similarly, the \emph{weak inverse linear dilatation function} of $f$ at $x$ is 
\begin{align*}
h_f^*(x)=\liminf_{s\to 0}H_f^*(x,s).
\end{align*}

\begin{defn}[Inverse metrically quasiregular mappings]\label{def:inverse metric quasiregular map}
	A branched covering $f\colon \tilde{\Omega}\to \Omega$ between two metric measure spaces is termed \textit{inverse metrically $H$-quasiregular} if the inverse linear dilatation function $H_f^*$ is finite everywhere and essentially bounded from above by $H$. 
\end{defn}

\begin{defn}[Inverse weak metrically quasiregular mappings]\label{def:inverse weak metric quasiregular map}
	A branched covering $f\colon \tilde{\Omega}\to \Omega$ between two metric measure spaces is termed \textit{inverse weak metrically $H$-quasiregular} if it satisfies 
	\begin{itemize}
		\item[i).] $h_f^*(x)<\infty$ for all $x\in \tilde{\Omega}$;
		\item[ii).] $h_f^*(x)\leq H$ for $\mu$-almost every $x\in \tilde{\Omega}$. 
	\end{itemize}
\end{defn}

Using the pullback factorization, we introduce the following new class of inverse analytically quasiregular mappings.

\begin{defn}[Inverse analytically quasiregular mappings]\label{def:inverse analytic def for qr}
	A branched covering $f\colon \tilde{\Omega}\to \Omega$ is said to be \textit{inverse analytically $K$-quasiregular} if $g^{-1}\in N^{1,Q}_{loc}(\tilde{\Omega}^f,\tilde{\Omega})$ and 
	\begin{align*}
	|\nabla g^{-1}|(z)^Q\leq KJ_{g^{-1}}(z)
	\end{align*}
	for $\lambda$-a.e. $z\in \tilde{\Omega}^f$. 
\end{defn}

The inverse geometric definition also relies on certain inequalities for the modulus of curve families.

\begin{defn}[Inverse geometric quasiregular mappings]\label{def:inverse geometric def for qr}
	A branched covering $f\colon \tilde{\Omega}\to \Omega$ is said to be \textit{inverse geometrically $K$-quasiregular} if it satisfies the \textit{$K_I$-inequailty} or the \textit{Poletsky's inequality}, i.e., for every curve family $\Gamma$ in $\tilde{\Omega}$, we have
	\begin{align*}
	\Modd_Q(f(\Gamma))\leq K\Modd_Q(\Gamma).
	\end{align*}
\end{defn}

We also introduce the following strong inverse geometrically quasiregular mappings, which can be viewd as a generalized \textit{V\"ais\"al\"a's inequality}.
\begin{defn}[Strong inverse geometrically quasiregular mappings]\label{def:strong inverse geometric def for qr}
	A branched covering $f\colon \tilde{\Omega}\to \Omega$ is said to be a \textit{strong inverse geometric $K$-quasiregular mapping} if it satisfies the following \textit{generalized V\"ais\"al\"a's inequality}: For each open subset $\tilde{\Omega}_0\subset \tilde{\Omega}$, each curve family $\Gamma$ in $\tilde{\Omega}_0$, $\Gamma'$ in $\Omega$, and for each $\gamma'\in \Gamma'$, there are curves $\gamma_1,\dots,\gamma_m\in \Gamma$ and subcurves $\gamma_1',\dots,\gamma_m'$ of $\gamma'$ such that for each $i=1,\dots,m$, $\gamma_i'=f(\gamma_i)$, and for almost every $s\in [0,l(\gamma')]$, $\gamma_i(s)=\gamma_j(s)$ if and only if $i=j$, then 
	\begin{align*}
	\Modd_Q(\Gamma')\leq \frac{K}{m}\Modd_Q(\Gamma).
	\end{align*}
\end{defn}

\subsection{Equivalences of definitions of quasiregularity}

By~\cite[Theorem 9.8]{hkst01} and~\cite[Theorem 1.1]{bkr07}, when $f\colon \tilde{\Omega}\to \Omega$ is a homeomorphism, and $\tilde{\Omega}$ and $\Omega$ have  $Q$-bounded geometry, the inverse (metric, weak metric, analytic, geometric) definitions for $f$ are, quantitatively, the forward (metric, weak metric, analytic geometric) definitions for $f^{-1}$. Moreover, in this case, each of these definitions is further equivalent to the local \textit{quasisymmetry}, quantitatively.

Using the pullback factorization, we have 
\begin{thm}[Theorem A, \cite{gw16}]\label{thm:main thm}
The following conclusions hold:
\begin{itemize}
	\item[i).] $f$ is analytically $K_O$-quasiregular if and only if it is geometrically $K_O$-quasiregular. Similarly, $f$ is inverse analytically $K_I$-quasiregular if and only if it is strong inverse geometrically $K_I$-quasiregular.
	
	\item[ii).] If $\tilde{\Omega}$ and $\Omega$ are Ahlfors $Q$-regular, and $Y$ has $c$-bounded turning, then either of the following two conditions
	\begin{itemize}
		\item[a).] $h_f(x)\leq h$ for all $x\in \tilde{\Omega}$;
		\item[b).] $h_f^*(x)\leq h$ for all $x\in \tilde{\Omega}$,
	\end{itemize}
	implies that $f$ is analytically $K_O$-quasiregular and inverse analytically $K_I$-quasiregular, with both constants $K_O$ and $K_I$ depending only on the constant of Ahlfors $Q$-regularity, and on $c$ and $h$.
	
	\item[iii).] If both $\tilde{\Omega}$ and $\Omega$ have $Q$-bounded geometry, then all of the metric, geometric and analytic definitions are quantitatively equivalent.
\end{itemize}
\end{thm}
Theorem~\ref{thm:main thm} generalizes the earlier results of~\cite{hkst01,bkr07,w12proc,w14} about quasiconformal mappings in a natural form to that of quasiregular mappings. Moreover, if $f$ satisfies the $K_I$-inequality with $\nu(f(\mathcal{B}_f))=0$, then $f$ satisfies the standard V\"ais\"al\"a's inequality with the same constant~\cite{gw16}.

\subsection{Branched quasisymmetric mappings}

We next define a proper subclass of metrically quasiregular mappings that carry similar global metric information, but less restrictive as those BDD mappings.

\begin{defn}[Branched quasisymmmetric mappings]\label{def:branched qs}
	Let $f\colon X\to Y$ be a branched covering. We say that $f$ is \textit{branched quasisymmetric} (BQS) if there exists a homeomorphism $\eta\colon [0,\infty)\to [0,\infty)$ such that
	\begin{equation}\label{eq:def of BQS}
	\frac{\diam f(E)}{\diam f(F)}\leq \eta\Big(\frac{\diam E}{\diam F} \Big)
	\end{equation}
	for all intersected 
	continua $E, F\subset X$.
	\end{defn}
	
\begin{rem}
If $f\colon X\to Y$ is a homeomorhpism between two metric spaces that have bounded turning and if $f$ satisfies the inequality~\eqref{eq:def of BQS}, then $f$ is quasisymmetric, quantitatively. See~\cite{gw16}.
\end{rem}	
		
%
%
%
	Simiar as quasiconformal mappings are locally quasisymmetric in spaces of bounded geometry, we have the following branched version.
\begin{thm}\label{thm:qr with bounded multiplicity is branched qs}
	Let $f\colon X\to Y$ be a weak metrically $H$-quasiregular mapping such that $N=N(f,X)<\infty$. Assume that both $X$ and $Y$ have locally $Q$-bounded geometry. Then $f$ is locally $\eta$-branched quasisymmetric, quantitatively, with $\eta$ depending only on $H$, $N$, and the data of $X$ and $Y$. 
\end{thm}
	
	\begin{rem}\label{rmk:on qr not branched qs}
		In Theorem~\ref{thm:qr with bounded multiplicity is branched qs}, the homeomorphism $\eta$ depends, quantitatively on the multiplicity $N$. In general, one can not get rid of this dependence from the theorem, as the simple analytic function $z\mapsto z^k$ indicated.
	\end{rem}
	
\section{Concluding remarks}	
All the preceeding works on the theory of quasiregular mapping or quasiconformal mappings require the spaces in question to have quantitative bounded geometry (i.e., both Ahlfors regularity and the Poincar\'e inequality with quantitative data), in order to get a rich theory. 

On the other hand, much of the theory has been extended to the general equiregular subRiemannian manifolds~\cite{gnw15,gl15} or even certain non-equiregular subRiemannian manifolds~\cite{a15}. Equiregular subRiemannian manifolds do not necdessarily have quantitative bounded geometry and non-equiregular subRiemannian manifolds are merely doubling (typically not Ahlfors regular).

The standard assumptions nowadays in analysis on metric spaces are doubling and Poincar\'e (i.e., the metric measure spaces are doubling and support abstract Poincar\'e inequalities). It is of great interest to know whether one can build the theory of quasiconformal/quasiregular mappings in such metric spaces. In particular, we would like to know whether the quantitative equivalence of definitions of quasiconformality/quasiregularity as in Theorem~\ref{thm:main thm} remains valid, and is further quantitatively equivalent to the local quasisymmetry/branched quasisymmetry as in Theorem~\ref{thm:qr with bounded multiplicity is branched qs}, in the more general setting. Note that the question is unknown even for mappings from $X$ to $\R^2$~\cite[Question 17.3]{r14}.

\subsection*{Acknowledgment}
Many thanks to wonderful event ``International Conference on Complex Analysis and Related Topics, The 14th Romanian-Finnish Seminar", where part of this work has been done.


\begin{thebibliography}{1}

\bibitem{a15}	
C. Ackermann, \textit{An approach to studying quasiconformal mappings on generalized Grushin planes}, Ann. Acad. Sci. Fenn. Math. 40 (2015), no. 1, 305-320. 	
	
\bibitem{bkr07}
Z. Balogh, P. Koskela and S. Rogovin, \textit{Absolute continuity of quasiconformal mappings on curves}, Geom. Funct. Anal. 17 (2007), no. 3, 645-664.		
	
	
\bibitem{bh04}
M. Bonk and J. Heinonen, \textit{Smooth quasiregular mappings with branching}, Publ. Math. Inst. Hautes\'etudes Sci. No. 100 (2004), 153-170.






\bibitem{c99}
J. Cheeger, Differentiability of Lipschitz functions on metric measure spaces, Geom. Funct. Anal. 9 (1999), no. 3, 428-517.

\bibitem{g62}
F.W. Gehring, \textit{Rings and quasiconformal mappings in space}, Trans. Amer. Math. Soc. 103 (1962), 353-393. 



\bibitem{gl15}
C.Y. Guo and T. Liimatainen, \textit{Equivalence of quasiregular mappings on subRiemannian manifolds via the Popp extension}, preprint 2016.

\bibitem{gnw15}
C.Y. Guo, S. Nicolussi Golo and M. Williams, \textit{Quasiregular mappings between subRimannian manifolds}, preprint 2015.

\bibitem{gw15}
C.Y. Guo and M. Williams, \textit{Porosity of the branch set of discrete open mappings with controlled linear dilatation}, preprint 2015.

\bibitem{gw16}
C.Y. Guo and M. Williams, \textit{Geometric function theory: the art of pullback factorization}, preprint 2016.

\bibitem{hk00}
P. Hajlasz and P. Koskela, \textit{Sobolev met Poincar\'e}, Mem. Amer. Math. Soc. 145 (2000), no. 688.

\bibitem{h02}
J. Heinonen, \textit{The branch set of a quasiregular mapping}, Proceedings of the International Congress of Mathematicians, Vol. II (Beijing, 2002), 691-700, Higher Ed. Press, Beijing, 2002.

\bibitem{hk11}
J. Heinonen and S. Keith, \textit{Flat forms, bi-Lipschitz parameterizations, and smoothability of manifolds}, Publ. Math. Inst. Hautes \'Etudes Sci. No. 113 (2011), 1-37. 

\bibitem{hk95}
J. Heinonen and P. Koskela, \textit{Definitions of quasiconformality}, Invent. Math. 120 (1995), 61-79.

\bibitem{hk98}
J. Heinonen and P. Koskela, \textit{Quasiconformal maps in metric spaces with controlled geometry}, Acta Math. 181 (1998), 1-61.

\bibitem{hkst01}
J. Heinonen, P. Koskela, N. Shanmugalingam and J.T. Tyson, \textit{Sobolev classes of Banach space-valued functions and quasiconformal mappings}, J. Anal. Math. 85 (2001), 87-139.

\bibitem{hkst15}
J. Heinonen, P. Koskela, N. Shanmugalingam and J.T. Tyson, \textit{Sobolev spaces on metric measure spaces: an approach based on upper gradients}, Cambridge Studies in Advanced Mathematics Series, 2015.

\bibitem{hr98}
J. Heinonen and S. Rickman, \textit{Quasiregular maps $\textbf{S}^3\to \textbf{S}^3$ with wild branch sets}, Topology 37 (1998), no. 1, 1-24. 

\bibitem{hr02}
J. Heinonen and S. Rickman, \textit{Geometric branched covers between generalized manifolds}, Duke Math. J. 113 (2002), no. 3, 465-529.


\bibitem{hs02}
J. Heinonen and D. Sullivan, \textit{On the locally branched Euclidean metric gauge}, Duke Math. J. 114 (2002), no. 1, 15-41.

\bibitem{m14}
G. Martin, \textit{The theory of quasiconformal mappings in higher dimensions}, I. Handbook of Teichm\"uller theory. Vol. IV, 619-677, IRMA Lect. Math. Theor. Phys., 19, Eur. Math. Soc., Z\"urich, 2014. 


\bibitem{mrv69}
O. Martio, S. Rickman and J. V\"ais\"al\"a, \textit{Definitions for quasiregular mappings}, Ann. Acad. Sci. Fenn. Ser. A I No. 448 (1969), 40 pp.

\bibitem{mrv70}
O. Martio, S. Rickman and J. V\"ais\"al\"a, \textit{Distortion and singularities of quasiregular mappings}, Ann. Acad. Sci. Fenn. Ser. A I No. 465 (1970), 13 pp.

\bibitem{mrv71}
O. Martio, S. Rickman and J. V\"ais\"al\"a, \textit{Topological and metric properties of quasiregular mappings}, Ann. Acad. Sci. Fenn. Ser. A I No. 488 (1971), 31 pp.


\bibitem{m73}
G.D. Mostow, \textit{Strong rigidity of locally symmetric spaces}, Annals of Mathematics Studies, No. 78. Princeton University Press, Princeton, N.J.; University of Tokyo Press, Tokyo, 1973. v+195 pp.

\bibitem{r14}
K. Rajala, \textit{Uniformization of two-dimensional metric surfaces}, Invent. Math., to appear.

\bibitem{r66}
Yu.G. Reshetnyak, \textit{Estimates of the modulus of continuity for certain mappings}, Sibirsk. Mat. Z., 7 (1966), 1106-1114; English transl. in Siberian Math. J., 7 (1966), 879-886.

\bibitem{re89}
Yu.G. Reshetnyak, \textit{Space mappings with bounded distortion}, Translations of Mathematical Monographs, 73. American Mathematical Society, Providence, RI, 1989.

\bibitem{r93}
S. Rickman, \textit{Quasiregular Mappings}, Ergeb. Math. Grenzgeb. (3) 26, Springer, Berlin, 1993.

\bibitem{s00}
N. Shanmugalingam, \textit{Newtonian spaces: an extension of Sobolev spaces to metric measure spaces}, Rev. Mat. Iberoamericana 16 (2000), no. 2, 243-279.

\bibitem{s75}
J. Sarvas, \textit{The Hausdorff dimension of the branch set of a quasiregular mapping}, Ann. Acad. Sci. Fenn. Ser. A I Math. 1 (1975), no. 2, 297-307.

\bibitem{t98} 
J. Tyson, \textit{Quasiconformality and quasisymmetry in metric measure spaces}, Ann.
Acad. Sci. Fenn. Ser. A I Math., 23 (1998), 525-548.

\bibitem{t00}
J. Tyson, \textit{Analytic properties of locally quasisymmetric mappings from Euclidean
	domains}, Indiana Univ. Math. J., 49 (2000), no. 3, 995-1016.

\bibitem{t01}
J. Tyson, \textit{Metric and geometric quasiconformality in Ahlfors regular Loewner spaces}, Conform. Geom.
Dyn. 5 (2001), 21-73 (electronic).

\bibitem{w12proc}
M. Williams, \textit{Geometric and analytic quasiconformality in metric measure spaces}, Proc. Amer. Math. Soc. 140 (2012), no. 4, 1251-1266. 

\bibitem{w14}
M. Williams, \textit{Dilatation, pointwise Lipschitz constants, and condition $N$ on curves}, Michigan Math. J. 63 (2014), no. 4, 687-700. 



\end{thebibliography}
\end{document}